\documentclass[12pt]{amsart}
\usepackage{amssymb}
\usepackage{amsmath}
\numberwithin{equation}{section}
\usepackage{graphicx}
\usepackage{epstopdf}

\usepackage{color}
\usepackage{latexsym}
\usepackage{amsfonts}

\newcommand{\sgn}{{\rm sgn}}

\renewcommand{\S}{\mathfrak S}
\newcommand{\s}{\sigma}

\newcommand\F{{\mathbb{F}}}

\DeclareFontFamily{OT1}{rsfs}{}
\DeclareFontShape{OT1}{rsfs}{m}{n}{ <-7> rsfs5 <7-10> rsfs7 <10->
rsfs10}{} \DeclareMathAlphabet{\mycal}{OT1}{rsfs}{m}{n}

\newcommand\beq{\begin{eqnarray*}}
\newcommand\eeq{\end{eqnarray*}}
\newcommand\ben{\begin{enumerate}}
\newcommand\een{\end{enumerate}}
\newcommand\bit{\begin{itemize}}
\newcommand\eit{\end{itemize}}

\newtheorem{thm}{Theorem}[section]

\newtheorem{prop}[thm]{Proposition}

\newtheorem{cor}[thm]{Corollary}
\newtheorem{defn}[thm]{Definition}

\newcommand\des{{\rm des}}
\newcommand\exc{{\rm exc}}
\newcommand\inv{{\rm inv}}
\newcommand\maj{{\rm maj}}
\newcommand\ai{{\rm ai}}
\newcommand\aid{{\rm aid}}
\newcommand\sg{{\mathfrak S}}
\newcommand\Des{{\rm DES}}
\newcommand\Exd{{\rm EXD}}
\newcommand\Exc{{\rm EXC}}

\newcommand\ch{{\rm ch}}

\newcommand\x{{\mathbf x}}
\newcommand\bbar{{\rm bar}}

\title[Excedance Number and Major Index]{$q$-Eulerian Polynomials: Excedance Number and Major index}
\author[Shareshian]{John Shareshian$^1$}
\address{Department of Mathematics, Washington University, St. Louis, MO}
\thanks{$^{1}$Supported in part by NSF Grants
 DMS 0300483 and DMS 0604233, and the Mittag-Leffler Institute}
\email{shareshi@math.wustl.edu}

\author[Wachs]{Michelle L. Wachs$^2$}
\address{Department of Mathematics, University of Miami, Coral Gables, FL 33124}
\email{wachs@math.miami.edu}
\thanks{$^{2}$Supported in part by NSF Grants
DMS 0302310 and DMS 0604562, and the Mittag-Leffler Institute}

\date{\today \\ \small MR Subject Classifications: 05A30, 05E05, 05E25}

\begin{document}

\maketitle
\begin{abstract} In this research announcement we present a new
$q$-analog of a classical formula for the exponential generating
function of the Eulerian polynomials.  The Eulerian polynomials
enumerate  permutations according to their number of descents or
their number of excedances.  Our $q$-Eulerian polynomials are the
enumerators for the joint distribution of the excedance statistic
and the major index.  There is a vast literature on $q$-Eulerian
polynomials which involve other combinations of  Mahonian and
Eulerian permutation statistics, but the combination of major
index and excedance number seems to have been completely
overlooked until now.   We use symmetric function theory to prove
our formula.  In particular, we prove a symmetric function version
of  our formula, which involves an intriguing new class of
symmetric functions.  We also present connections with
representations of the symmetric group on the homology of a  poset
recently introduced  by Bj\"orner and Welker and on the cohomology
of the toric variety associated with the Coxeter complex of the
symmetric group, studied by  Procesi, Stanley,  Stembridge,
Dolgachev and Lunts. \end{abstract}

\section{Introduction} The subject of permutation statistics
originated in the  early 20th century work of Major Percy
MacMahon \cite{mac1,mac2} and has developed into an active and
important area of enumerative combinatorics over the last four
decades.   It deals with the enumeration of  permutations
according to   natural statistics.    A permutation statistic is
simply a function from the symmetric group $\S_n$ to the set of
nonnegative integers.     MacMahon studied four fundamental
permutation statistics,  the inversion index, the major index, the
descent number and the excedance number, which we define below.

Let $[n]$ denote the set $\{1,2,\dots,n\}$.   For each $\sigma \in
\S_n$, the descent set of $\sigma$ is defined to be
$$\Des(\sigma):= \{i \in [n-1] : \sigma(i) > \sigma(i+1) \},$$ and
the excedance set is defined to be $$\Exc(\sigma):= \{i \in [n-1]
: \sigma(i) > i\}.$$   The  descent number and excedance number
are defined respectively by $$ \des(\sigma) := |\Des(\sigma)|
\qquad \mbox{ and }\qquad  \exc(\sigma) := |\Exc(\sigma)|.$$ For
example, if $\sigma = 32541$, written in one line notation, then
$$\Des(\s) = \{1,3,4\} \qquad \mbox{ and }\qquad \Exc(\s) = \{1,3\};$$
 hence $\des(\s) = 3 $ and $\exc(\s) = 2$.  If $i \in \Des(\s)$ we say
 that $\s$ has a descent at $i$.  If $i \in \Exc(\s)$ we say that
 $\s(i)$ is an excedance of $\s$.

MacMahon \cite{mac1,mac2} observed that the descent number and
excedance number are equidistributed, that is, the number of
permutations in $\S_n$ with $j$ descents equals the number of
permutations with $j$ excedances for all $j$.   (There is a well-known
combinatorial proof of this fact due to Foata \cite{foa1,foa4}.)  These
numbers were first studied by Euler   and have come to be known  as
the Eulerian numbers.  They are the coefficients of the {\em
Eulerian polynomials} $$A_n(t) := \sum_{\sigma \in \sg_n}
t^{\des(\sigma)} =  \sum_{\sigma \in \sg_n} t^{\exc(\sigma)}. $$
Any permutation statistic that is equidistributed with $\des$ and
$\exc$ is said to be an {\em Eulerian statistic}.

The Eulerian numbers and the Eulerian polynomials have been
extensively studied in many different contexts in the mathematics
and computer science literature.  For excellent treatments of this subject, see the classic lecture notes of Foata and Schutzenberger \cite{fs}, 
the recent lecture notes of Foata and Han \cite{fh},  and Section~5.1 of Knuth's classic book series ``The Art of Computer Programming'' \cite{kn}. 
The exponential generating function formula,
\begin{equation} \label{expgen}\sum_{n\ge 0} A_n(t) {z^n \over n!}
= {1-t \over e^{z(t-1)} -t} \end{equation}
where $A_0(t) =1$, is attributed to Euler in \cite{kn}.

The major index of a permutation is defined  by
$$\maj(\sigma) := \sum_{i\in \Des(\sigma)} i .$$
 MacMahon \cite{mac2} proved
that the major index is equidistributed with the inversion statistic
$$ \inv(\sigma) := |\{(i,j) : 1 \le i<j\le n \,\,\&
\,\, \sigma(i) > \sigma(j) \}|$$ and Rodrigues \cite{rod} proved
the second equality in \begin{equation} \label{mah} \sum_{\sigma
\in \sg_n} q^{\maj(\sigma)} = \sum_{\sigma \in \sg_n}
q^{\inv(\sigma)} = [n]_q!,\end{equation} where $$[n]_q := 1+q
+\dots +q^{n-1}$$ and $$[n]_q! :=  [n]_q [n-1]_q \cdots [1]_q.$$
(An elegant combinatorial proof of the first equality  in
(\ref{mah}) was obtained by Foata \cite{foa2,foa4}.)
 Any permutation statistic  that is equidistributed with the major
index and inversion index is said to be a {\em Mahonian
statistic}.

Note that by setting $q=1$ in (\ref{mah}), one gets the formula
$n!$ for the number of permutations. Equation (\ref{mah}) is a
beautiful ``$q$-analog'' of this formula and is the fundamental example 
of the subject of permutation statistics and
$q$-analogs, in which one seeks to obtain nice $q$-analogs of
 enumeration formulas.

One can look for nice $q$-analogs of the Eulerian polynomials by
considering the joint distributions of the Mahonian and Eulerian
statistics given above.   Consider the four possibilities,
$$A_n^{\inv,\des}(q,t) := \sum_{\sigma \in \sg_n} q^{\inv(\sigma)} t^{\des(\sigma)} $$
$$A_n^{\maj,\des}(q,t) := \sum_{\sigma \in \sg_n} q^{\maj(\sigma)}  t^{\des(\sigma)} $$
$$A_n^{\inv,\exc}(q,t) := \sum_{\sigma \in \sg_n}q^{\inv(\sigma)} t^{\exc(\sigma)} $$
$$A_n^{\maj,\exc}(q,t) := \sum_{\sigma \in \sg_n} q^{\maj(\sigma)} t^{\exc(\sigma)} .$$
There are many interesting results  on the first three
$q$-Eulerian polynomials and on multivariate distributions of all
sorts of combinations of Eulerian and Mahonian statistics (for a sample see
\cite{br,csz, foa3,fs2,fz, gg, rrw, ra, sk,st1,wa2}).  These include   Stanley's \cite{st1}  $q$-analog
of (\ref{expgen}) given by,
$$\sum_{n\ge 0} A_n^{{\inv,\des}}(q,t) {z^n \over [n]_q!} = {(1-t) \over \exp_q(z(t-1)) -t}$$
where $$\exp_q(z) := \sum_{n\ge 0}  {z^n \over [n]_q!}.$$

Surprisingly, we have found no mention of the fourth $q$-Eulerian
polynomial $A_n^{\maj,\exc}(q,t)$ anywhere in the literature. Here
we announce the following remarkable $q$-analog of
(\ref{expgen}).

\begin{thm}\label{expgenth} The $q$-exponential
generating function for
$A_n^{\maj,\exc}(q,t)$ is given by
\begin{equation}  \label{expgeneq} \sum_{n\ge 0} A_n^{\maj,\exc}(q,t) {z^n \over [n]_q!}  = {(1-tq)
\exp_q(z) \over \exp_q(ztq) - tq \exp_q(z)}\,\, , \end{equation} where
$A_0^{\maj,\exc}(q,t) = 1$.
\end{thm}

When $q=1$, the formula (\ref{expgeneq})   reduces to
(\ref{expgen}) since $${(1-tq) \exp_q(z)  \over \exp_q(ztq) -tq
\exp_q(z)} = {(1-t) e^z \over e^{zt} -te^z} = {(1-t) \over
e^{z(t-1)} -t}.$$ Though not quite as easily,  one can show that
when   $t=1$, the formula (\ref{expgeneq}) reduces to (\ref{mah}).

In the process of proving Theorem~\ref{expgenth}, we obtained the
following result.

\begin{thm} \label{stem} Let  ${\rm fix}(\sigma)$ denote the number
of fixed points of $\s\in \S_n$, i.e., the number of $i \in [n]$
such that $\s(i) = i$.  Then
$$\sum_{\sigma
\in
\sg_n} q^{\maj(\sigma)}t^{\exc(\sigma)}r^{{\rm fix}(\sigma)}  = $$
$$  \sum_{m = 0}^{\lfloor {n \over 2} \rfloor} (tq)^m \!\!\!\!\sum_{\scriptsize
\begin{array}{c} k_0\ge 0  \\ k_1,\dots, k_m \ge 2 \\ \sum k_i = n
\end{array}} \left[\begin{array}{c} n \\k_0,\dots,k_m\end{array}\right]_q\,\,
r^{k_0}
\prod_{i=1}^m [k_i-1]_{tq},$$
where $$\left[\begin{array}{c} n \\k_0,\dots,k_m\end{array}\right]_q = {[n]_q!
\over [k_0]_q![k_1]_q!\cdots [k_m]_q!}.$$
\end{thm}

In the next section we describe the techniques that we used to
prove these theorems.  They involve an interesting class of
symmetric functions and a symmetric function identity which
generalizes  Theorem~\ref{expgenth}. We prove this identity by
devising an interesting analog of a necklace construction of
Gessel and Reutenauer \cite {gr} and by generalizing a bijection
of Stembridge \cite{stem1}.

In  Section~\ref{rep} we discuss a connection with two  graded
representations of the symmetric group, which turn out to be
isomorphic.  We show that a specialization of the Frobenius
characteristic of these representations yields
$A^{\maj,\exc}(q,t)$.  One of the representations is the
representation of the symmetric group on the cohomology of the
toric variety associated with the Coxeter complex of the symmetric
group.  This representation  was studied by Procesi \cite{pr},
Stanley \cite{st2}, Stembridge \cite{stem1}, \cite{stem2}, and
Dolgachev and Lunts \cite{dl}. The other representation is the
representation of the symmetric group on the homology of maximal
intervals of a certain intriguing poset introduced by Bj\"orner
and Welker \cite{bw} in their study of connections between poset
topology and commutative algebra.  In fact,  our study of the
latter representation is what led us to discover   formula
(\ref{expgenth}) and its symmetric function generalization, in the
first place.

Various authors have studied Mahonian (resp. Eulerian) partners to
Eulerian (resp. Mahonian) statistics whose joint distribution is
equal to a known Euler-Mahonian distribution.   We mention, for
example, Foata \cite{foa3},  Foata and Zeilberger \cite{fz},
Skandera \cite{sk}, and Clarke, Steingrimisson and Zeng
\cite{csz}.   In  Section~\ref{new} we define a new Mahonian
statistic to serve as a partner for $\des$ in the $(\maj, \exc)$
distribution.   We do not have a simple proof of the
equidistribution.  We have a highly nontrivial proof which uses
tools  from poset topology  and the symmetric function results
announced in Sections~\ref{symsec} and \ref{rep}.

Details of the proofs discussed in this announcement, as well as further consequences and open problems, will appear in a forthcoming paper.

\section{Symmetric function generalization} \label{symsec}
In this section we present a symmetric function generalization of
Theorem~\ref{expgenth}.

Let
$$H(z) = H(\mathbf x, z) := \sum_{n\ge 0} h_n(\x) z^n,$$
where $h_n(\mathbf x)$ denotes the complete homogeneous symmetric
function in the indeterminates $\mathbf x = (x_1,x_2,\dots)$, that
is $$h_n(\mathbf x) := \sum_{1 \le i_1 \le i_2 \le \dots \le i_n}
x_{i_1}x_{i_2} \dots x_{i_n}$$ for $n \ge 1$, and $h_0 =1$. By
setting $x_i := q^{i-1}$, for all $i$, and $z := z(1-q)$ in
$H(\x,z)$, one obtains $\exp_q(z)$, see \cite{st3}. It follows
that  \begin{equation}\label{Heq}  \left .{(1-t) H(\x,z) \over
H(\x,zt) -tH(\x,z)} \right |_{\scriptsize \begin{array}{l}
x_i:=q^{i-1} \\ z:= z(1-q) \end{array}} = {(1-t) \exp_q(z) \over
\exp_q(zt) -t\exp_q(z)} .\end{equation}

We will construct for each $n,j\ge 0$, a quasisymmetric function $Q_{n,j}(\x)$ whose
generating function $\sum_{n,j \ge 0} Q_{n,j}(\x) t^j z^n$
specializes to
$$\sum_{n\ge 0} \sum_{\s \in \S_n} q^{\maj(\s) - \exc(\s)} t^{\exc(\s)} {z^n \over [n]_q!}$$
when we set $x_i := q^{i-1}$ and $z := z(1-q)$.  Thus  by taking
specializations of both sides of (\ref{symgen}) below and setting
$t:=tq$, we obtain (\ref{expgeneq}).

For $\s \in \sg_n$, let $\bar \s$ be the barred word obtained from
$\s$ by placing a bar above each {excedance}. For example, if $\s
= 531462$ then {$\bar \s = \bar 5 \bar 3 14 \bar 6 2$.} View $\bar
\s$ as a word over ordered alphabet $$\{ \bar 1 < \bar 2 < \dots <
\bar n < 1 <  2 < \dots < n\}.$$ We extend the definition of
descent set   from permutations to words $w $ of length $n$ over
an ordered alphabet by letting
$$\Des(w) := \{i \in [n-1] : w_i > w_{i+1}\},$$
where $w_i$ is the $i$th letter of $w$.  Now define the
excedance-descent set of a permutation $\s \in \S_n$ to be
$$\Exd(\s) := \Des(\bar \s).$$
For example,
$\Exd({531462}) = \Des({\bar 5} \bar 3 14{\bar 6 2}) = \{1,4\}$.
The interesting thing about $\Exd$ is that  for all $\s \in \S_n$,
\begin{equation}\label{exd}  \sum_{i \in \Exd(\s)} i = \maj(\s) - \exc(\s).\end{equation}

For $S \subseteq [n-1]$ and $n \ge 1$, define the quasisymmetric
function $$F_{S,n}(x_1,x_2,\dots):= \sum_{\scriptsize
\begin{array}{c} i_1 \ge \dots \ge i_n \\ j \in S \Rightarrow i_j
> i_{j+1} \end{array}} x_{i_1} \dots x_{i_n},$$  and let
$F_{\emptyset, 0} = 1$.

A basic result in  Gessel's theory of  quasisymmetric functions
(see eg., \cite{st3}) is that
$$F_{S,n}(1,q,q^2,\dots) = {q^{\sum_{s \in S} s} \over (1-q) (1-q^2) \dots (1-q^n)}.$$
Hence  it follows from (\ref{exd}) that for  all $\s \in \S_n$,
$$F_{\Exd(\s),n} (1,q,q^2,\dots) = {q^{\maj(\s) -\exc(\s) } \over (1-q) (1-q^2) \dots (1-q^n)}.$$
For any $n,j \ge 0$, let $$ Q_{n,j} = Q_{n,j}(\x):=
\sum_{\scriptsize \begin{array}{c} \s \in \sg_n \\ \exc(\s) = j
\end{array}} F_{\Exd(\s),n}(\x).$$   By taking the specialization
of the generating function we get,
\begin{equation}
\label{Qeq}\left . \sum_{n,j \ge 0} Q_{n,j}(\x) t^j z^n \right
|_{\scriptsize \begin{array}{l} x_i:=q^{i-1} \\ z:= z(1-q)
\end{array}} = \sum_{n\ge 0} \sum_{\s \in \S_n} q^{\maj(\s) -
\exc(\s)} t^{\exc(\s)} {z^n \over [n]_q!} .\end{equation}

It follows from (\ref{Heq}) and (\ref{Qeq}) that by setting
$x_i:=q^{i-1} , z:= z(1-q)$ and $t:=tq$ in the following result we
obtain  Theorem~\ref{expgenth}.
\begin{thm} \label{symgenth}
\begin{equation} \label{symgen}
\sum_{n,j \ge 0} Q_{n,j} t^j z^n = {(1-t) H(z) \over H(zt) -tH(z)}.
\end{equation}
\end{thm}

The proof of this theorem requires an alternative characterization
of $Q_{n,j}$ which involves an interesting analog of  a
construction of  Gessel and Reutenauer \cite{gr}.   Gessel and
Reutenauer deal with circular words over the alphabet of positive
integers.  We  consider circular words over the alphabet  of
barred and unbarred positive integers.   For each such circular
word and any starting position, one gets an infinite word by
reading the circular word in a clockwise direction. If one gets a
distinct infinite word for each starting position, then the
circular word is said to be {\em primitive}.  For example $(\bar
1,1,1)$ is primitive while $(\bar 1, 2, \bar 1, 2)$ is not. The
{\em absolute value} of a letter is the letter obtained by erasing
the bar if there is one. We will say that a  primitive  circular
word is  a {\em necklace} if each letter that is followed
(clockwise) by a letter greater in absolute value is barred   and
each letter that is followed by a letter smaller in absolute value
is unbarred. Letters that are followed by letters equal in
absolute value have the option of being barred or not.  A circular
word  consisting of one barred letter is not a necklace.  For example the following
circular words are necklaces: $$(\bar 1 , 3, 1, \bar 1, 2, 2),
(\bar 1 , 3, \bar 1, \bar 1, 2, 2),  (\bar 1 , 3, 1, \bar 1, \bar
2,  2), (\bar 1 , 3, \bar 1, \bar 1, \bar 2, 2), (3),$$ while
$(\bar 1 , \bar 3, 1,  1, 2, \bar 2)$ and $(\bar 3)$ are not.

Again we will need to order the barred letters, but this time  by
$$\bar 1 <1 <  \bar 2 < 2 < \dots .$$
We order the necklaces by lexicographic order of the
lexicographically smallest infinite word obtained by reading the
necklace in a clockwise direction at some starting position.  An
{\em ornament}  is a weakly decreasing finite sequence of
necklaces.  The type $\lambda(R)$ of an ornament $R$ is the
partition whose parts are the sizes of the necklaces  in $R$. The
weight $w(R)$ of an ornament $R$  is the product of the weights of
the letters of  $R$, where the weight of the letter $a$ is  the
indeterminate $x_{|a|}$, where $|a|$ denotes the absolute value of
$a$.  For example
$$\lambda(  (\bar 1, 2, 2),(\bar 1, \bar 2, 3,3,2)) = (5,3)$$ and
$$w(  (\bar 1, 2, 2),(\bar 1, \bar 2, 3,3,2)) = x_1^2 x_2^4 x_3^2.$$
For each partition $\lambda$ and nonnegative integer $j$, let
$\mathfrak R_{\lambda,j}$ be the set of ornaments of type
$\lambda$ with $j$ bars.

\begin{thm} \label{ornth}
For all $\lambda\vdash n$ and $j= 0,1,\dots,n-1$, let
$$Q_{\lambda,j} =  \sum_{\s} F_{\Exd(\s),n}$$ summed over all
permutations of cycle type $\lambda$ with $j$ excedances. Then
$$Q_{\lambda,j} = \sum_{R \in \mathfrak R_{\lambda,j}} w(R).$$
\end{thm}

This theorem is proved via a bijection between ornaments of type
$\lambda$ and    permutations of cycle type $\lambda$ paired with
``compatible''  weakly increasing sequences of positive integers.
The theorem has several interesting consequences.  For one thing,
it can be  used it to prove that the quasisymmetric functions
$Q_{\lambda,j}$ and $Q_{n,j}$ are actually symmetric.   It also
has the following useful consequence.

\begin{cor}
\label{dercor} Let  $$\tilde Q_{n,j} =
\sum_{\scriptsize\begin{array}{c} \s \in \mathcal D_n\\ \exc(\s) =
j  \end{array}} F_{\Exd(\s),n}$$ where $\mathcal D_n$ is the set
of derangements in $\S_n$. Then $$ Q_{n,j} = \sum_{k = 0}^n h_k
\tilde Q_{n-k,j}.$$
\end{cor}

It follows from Corollary~\ref{dercor} that Theorem~\ref{symgen}
is equivalent to
\begin{equation} \label{symgen2}
\sum_{n,j \ge 0} \tilde Q_{n,j} t^j z^n = {1-t  \over H(zt) -tH(z)},
\end{equation}
which in turn, is equivalent to the recurrence relation
\begin{equation} \label{rr}
\tilde Q_{n,j} = \sum_{\scriptsize \begin{array}{c}0 \le m \le n-2
\\ j+m-n < i < j \end{array}} \tilde Q_{m,i}   h_{n-m}.
\end{equation}
We establish this recurrence relation by introducing another type
of configuration, closely related to ornaments.

Define a {\em banner}  $B$ to be a word over the alphabet of
barred and unbarred positive integers, where $B(i) $ is barred if
$|B(i)| < |B(i+1)|$ and $B(i)$ is unbarred if $|B(i)| > |B(i+1)|$
or $i = length(B)$. All other letters have the option of being
barred.   The weight of a banner is the product of the weights of
its letters.

A {\em Lyndon word} over an ordered alphabet is a  word that is
strictly lexicographically smaller than all its circular rearrangements.  A
{\em Lyndon factorization} of a word over an ordered alphabet is a
factorization into a weakly lexicographically  decreasing sequence
of  Lyndon words.  It is a result of Lyndon \cite{l} that every
word has a unique Lyndon factorization.  The Lyndon type of a word
is the partition whose parts are the lengths of the words in its
Lyndon factorization.  For each partition $\lambda$ and positive
integer $j$, let $\mathfrak B_{\lambda,j}$ be the set of banners
with $j$ bars whose Lyndon type is $\lambda$.

By turning the Lyndon words in the Lyndon factorization of a
banner into circular words, we obtain an ornament.  This map from
banners to ornaments is  the  bijection whose existence is
asserted in  the following proposition.

\begin{prop} For any partition $\lambda$ and nonnegative integer $j$, there is a weight-preserving bijection from
$\mathfrak B_{\lambda,j}$ to $\mathfrak R_{\lambda,j}$. \end{prop}

\begin{cor} \label{ban} Let $\tilde {\mathfrak B}_{n,j}$ be the
set of banners of length $n$ with $j$ bars whose  Lyndon type has
no parts of size $1$.  Then
 $$\tilde Q_{n,j} = \sum_{B \in \tilde {\mathfrak B}_{n,j} }w(B) $$
\end{cor}

Define a {\em marked sequence} $(\alpha,j)$  to be a weakly
increasing finite sequence $\alpha$ of positive integers together
with an integer  $j$ such that $1 \le j \le \mbox
{length}(\alpha)-1$.  Let $\mathfrak M_n$ be the set of marked
sequences of length $n$ and let $\tilde{\mathfrak B}_n$ be the set
of banners of length $n$ whose Lyndon type has no parts of size
$1$.

\begin{thm} \label{bij} For all $n > 0$, there is a bijection
$$\gamma: \tilde {\mathfrak B}_{n} \to \bigcup_{m>0}
\tilde{\mathfrak B}_{m} \times \mathfrak M_{n-m},$$ such that  if
$\gamma(B) = (B^\prime,(\alpha,j)) $ then $$w(B) =
w(B^\prime)w(\alpha)$$ and $$\bbar(B) = \bbar(B^\prime) + j,$$
where $\bbar(B)$ denotes the number of bars of $B$. \end{thm}

We will not describe the bijection here except to say that, when
restricted to banners with distinct letters, it reduces to   a
bijection from permutations to marked words that Stembridge
\cite{stem1} constructed to study the representation of the
symmetric group on the cohomology of the toric variety assoiciated
with the type A Coxeter complex.  (We discuss this representation
in Section~\ref{rep}.)  Banners in  $ \tilde {\mathfrak B}_{n}$
admit a certain kind of decomposition, called a decreasing
decomposition in \cite{dw}. The decreasing decomposition plays the
role in our bijection that the cycle decomposition of permutations
plays in Stembridge's bijection. Corollary~\ref{ban} and
Theorem~\ref{bij} are all that is needed to establish the
recurrence relation (\ref{rr}), which yields our main result,
Theorem~\ref{symgenth}.

 \section{Some Representation Theoretic Consequences} \label{rep}

The Frobenius characteristic $\ch$ is a fundamental homomorphism
from the ring of  representations of symmetric groups to the ring
of symmetric functions.  In this section we present two
representations  whose Frobenius characteristic is   $Q_{n,j}$.

The first representation involves  the toric variety associated
with the   Coxeter complex of a Weyl group.   
Let $X_n$ be  the
toric variety associated with the Coxeter complex of $\sg_n$. The
action of $\sg_n$ on $X_n$ induces a representation of $\sg_n$ on
the cohomology $H^{2j}(X_n)$ for each $j = 0,\dots,n-1$. (Cohomology in odd degree vanishes.)
Stanley \cite{st2}, using a formula of Processi \cite{pr}, proves
that
$$\sum_{n\ge 0} \sum_{j=0}^{n-1} \ch H^{2j}(X_n)\,t^{j} z^n
= {(1-t) H(z) \over H(zt) -tH(z)}. $$  Combining this with
Theorem~\ref{symgenth} yields the following conclusion.

\begin{thm} For all $j = 0,1, \dots, n-1$,
$$\ch H^{2j}(X_n)=Q_{n,j}.$$
\end{thm}

The second representation involves poset topology, a subject  in
which topological properties of a simplicial complex associated
with a poset are studied, see  \cite{w1}.  The faces of the
simplicial complex, called the order complex of the poset, are the
chains of the poset.   Here we consider the homology of the order
complex of the Rees product of two simple posets. The Rees product
is a poset construction recently introduced by Bj\"orner and
Welker \cite{bw} in their study of relations between poset
topology and commutative algebra.

\begin{defn}   Let $P$ and $Q$ be pure (ranked) posets with respective rank
functions $r_P$ and $r_Q$.  The Rees product $P*Q$ of $P$ and $Q$ is defined as
follows:
$$P*Q :=\{(p,q) \in P \times Q : r_P(p) \ge r_Q(q)\}$$
with order relation  given by $(p_1,q_1) \le (p_2,q_2)$ if the
following holds \bit \item $ p_1 \le_P p_2 $ \item $q_1 \le_Q q_2$
\item $r_P(p_2) -r_P(p_1) \ge r_Q(q_2) -r_Q(q_1)$ \eit
\end{defn}

Let $B_n$ be the Boolean algebra (ie., the lattice of subsets of
$[n]$ ordered by inclusion) and let $C_n$ be the chain $1<2< \dots
< n$.   The maximum elements of $(B_n \setminus \{\emptyset\}) *
C_n $ are of the form $([n],j)$, where $j=1,\dots, n$.  Let
$I_{n,j}$ be the set of elements of  $(B_n \setminus
\{\emptyset\}) * C_n $ that are smaller than $([n],j)$ and  let
$\tilde H_{i}(I_{n,j})$ be the reduced simplicial (complex)
homology of the order complex of  $I_{n,j}$.  It follows from
results of Bj\"orner and Welker that homology vanishes below the
top dimension $n-2$.  The symmetric group $\sg_n$ acts on $I_{n,j}$  in an obvious way and this induces a representation on $\tilde H_{n-2}( I_{n,j})$.  We prove the following result using
techniques from poset topology.

\begin{thm} \label{rees}
$$1+ \sum_{n\ge 1} \sum_{j=1}^n \ch (\tilde H_{n-2}(I_{n,j})
\otimes \sgn) \,t^{j-1} z^n = {(1-t) H(z) \over H(zt) -tH(z)},$$
where $\sgn$ denotes the sign representation. Consequently  for
all $n,j$, $$\ch (\tilde H_{n-2}( I_{n,j})\otimes \sgn) =
Q_{n,j-1}$$ and as $\S_n$-modules $$\tilde H_{n-2}(
I_{n,j})\otimes \sgn \cong H^{2j}(X_n).$$
\end{thm}

We conjecture that for all $\lambda$ and $j$, the symmetric
function $Q_{\lambda,j}$ is also the Frobenius characteristic of
some representation.  One consequence of Theorem~\ref{ornth} is
that $Q_{\lambda,j}$ can be described as a product of  plethysms
of symmetric functions of the form $Q_{(n),i}$, where $(n)$
denotes a partition with a single part.  Hence if the conjecture
holds for all $Q_{(n),i}$ then it holds in general. We  use
ornaments and banners to show that if the conjecture does hold
then the restriction to $\S_{n-1}$ of the representation whose
Frobenius characteristic is $Q_{(n),i}$,  has Frobenius
characteristic  $Q_{n-1,i-1}$.

\section{A new Mahonian statistic} \label{new}

In this section we describe a new Mahonian statistic whose joint
distribution with $\des$ is the same as the joint distribution of
$\maj$ and $\exc$.

 An {\em admissible inversion} of $\s \in \sg_n$ is a pair
 $(\s(i),\s(j))$ such that the following conditions hold:
\bit
\item $i <j$
\item $\s(i) > \s(j)$
\item either
    \bit
        \item[$\circ$]  $\s(j) < \s(j+1)$ or
        \item[$\circ$]  $\exists k $ such that $i<k<j$ and $\s(k) < \s(j)$.
        \eit
\eit
Let $\ai(\s):= \#$ admissible inversions of $\s$.
Define the statistic $${\aid}(\s) := \ai(\s) + \des(\s).$$
For example, the
admissible inversions of $24153$ are $(2,1), (4,1)$ and $(4,3).$
So  $\aid(24153)= 3+2$.

\begin{thm} \label{aid} For all $n \ge 1$,
$$\sum_{\s \in \S_n} q^{\aid(\s)} t^{\des(\s)}=
\sum_{\s \in \S_n} q^{\maj(\s)} t^{\exc(\s)}.$$
\end{thm}

We do not have a direct proof of this simple identity except when $t$ or $q$ is 1.  Our proof
relies on Theorem~\ref{expgenth},  a $q$-analog of
Theorem~\ref{rees}, and techniques from poset topology.   We
consider the Rees product $(B_n(q) \setminus \{(0)\}) * C_n$,
where $B_n(q)$ is the lattice of subspaces of the vector space
$\F_q^n$.  Let $I_{n,j}(q)$ be the set of elements in $(B_n(q)
\setminus \{(0)\}) * C_n $ that are less than the maximal element
$(\F_q^n,j)$.  We first use a well-known tool from poset topology,
called lexicographic shellability, to prove that
\begin{equation}\label{lex} \dim \tilde H_{n-2}(I_{n,j}(q))
={\sum_{\scriptsize{\begin{array}{c} \sigma \in \sg_n
\\
\des(\sigma) = j-1\end{array}}} \hspace{-.2in} q^{\ai(\s)}}.
\end{equation}
We then use other tools from poset topology to prove a theorem
analogous to Theorem~\ref{rees} which states that
\begin{equation}
\sum_{n\ge 0} \sum_{j=1}^n\dim \tilde H_{n-2}(I_{n,j}(q))
t^{j-1}{z^n \over [n]_q!}  = {(1-t) \exp_q(z) \over \exp_q(zt) - t
\exp_q(z)} .\end{equation} Theorem~\ref{aid} now follows from
Theorem~\ref{expgenth} and equation (\ref{lex}).

\section{Acknowledgements}
The research presented here began while both authors were visiting the Mittag-Leffler Institute as participants in a combinatorics program organized by Anders Bj\"orner and Richard Stanley.  We thank the Institute for its hospitality and support.  We are also grateful to Ira Gessel for some very useful
discussions.

\end{document}